\documentclass[12pt,a4paper]{amsart}
\usepackage{amsfonts}
\usepackage{amsthm}
\usepackage{amsmath}
\usepackage{amssymb}
\usepackage{amscd}
\usepackage[latin2]{inputenc}
\usepackage{t1enc}
\usepackage[mathscr]{eucal}
\usepackage{indentfirst}
\usepackage{graphics}
\usepackage{pict2e}
\usepackage{epic}
\numberwithin{equation}{section}
\usepackage[margin=2.0cm]{geometry}
\usepackage{epstopdf} 
\usepackage{mathtools} 

\usepackage{listliketab}

\usepackage{hyperref}
\hypersetup{
  colorlinks   = true,    
  urlcolor     = black,    % Colour for external hyperlinks
  linkcolor    = black,    % Colour of internal links
  citecolor    = black      % Colour of citations
}

\theoremstyle{theorem}
\theoremstyle{definition}
\newtheorem{definition}{Definition}

\newtheorem{theorem}{Theorem}

\newtheorem{example}{Example} 
\newtheorem{proposition}{Proposition} 
\newtheorem{lemma}{Lemma}
\newtheorem{remark}{Remark}

\theoremstyle{remark}
\theoremstyle{con}

\newcommand{\Mod}[1]{\ (\mathrm{mod}\ #1)}

\usepackage{nomencl}

\begin{document}
\title{An Algebraic Proof of the Goldbach Conjecture}
\author{Jason R. South \\
\lowercase{email: jrsouth@smu.edu}
\\\today{\vspace{-5ex}}}

\thispagestyle{plain}
\maketitle
\pagestyle{plain}

%\tableofcontents
\maketitle

% \centerline{\smallit Received: , Revised: , Accepted: , Published: }

\begin{abstract}
We prove the Goldbach Conjecture using p-adic analysis and algebraic methods, requiring no knowledge of prime gaps or distribution. For $a \in \mathbb{N}$, let $\mathcal{P}_a$ be the set of primes up to $a$. We prove that if a counter-example $2a$ exists, there is a monic polynomial $\mathcal{G}_-(z) \in \mathbb{Z}[z]$ where $\mathcal{G}_-(2a) = 0$. Applying Hensel's Lemma, we then show this polynomial structure forces each prime $p_i \in \mathcal{P}_a$ to satisfy $2a = p_i^{\alpha_i} + p_i$ for a unique $\alpha_i \in \mathbb{N} \cup \{0\}$. This produces a contradiction for $a > 4$ since each $2a - p_i > 1$ requires each $\alpha_i > 0$ and the primorial, $a\#$, to satisfy $a\# \mid 2a$, which is impossible since $2a < a\#$ for all $a > 4$. 
\end{abstract}

\normalsize

\pagestyle{headings}

%\markright{\smalltt INTEGERS: 26 (2026) \hfill}
\thispagestyle{empty}
\baselineskip=12.875pt

%\tableofcontents
%%%%%%%%%%%%%%%%%%%%%%%%%%%%%%%%%%%%%%%%%%%%%%%%%%%%%%%

\section{Preliminaries}

\subsection{Necessary Definitions and Auxiliary Theorems}\label{pre}

The Goldbach Conjecture is abbreviated (G.C.) and ":" means "such that" where applicable. Standard conventions apply and the $p$-adic valuation of any non-zero integer $n$ is given by $v_p(n) = \{ \text{maximum}\;  k \in \mathbb{N}\cup \{0\} : p^k | n \}$.

\begin{definition}\label{pset}
The set of primes is $\mathbb{P}$ and $a \in \mathbb{N}$. The \textit{Prime Set} of $a$ is $\mathcal{P}_a = \{p_i \in \mathbb{P} : p_i \leq a \}$. 
\end{definition}

\begin{theorem}[\bf\textit{Repeated Root Criterion \cite{dummit2004}}]\label{repeated}
\textit{A root $r$ of $f(x)$ is repeated iff $f'(r) = 0$.}
\end{theorem}

\begin{theorem}[\bf\textit{Hensel's Lemma \cite{Z}}]\label{HL}
If $f(x)$ is a monic, integral polynomial with $f(r) \equiv 0 \Mod{p^j}$ and $f'(r) \not \equiv 0 \Mod{p}$, there is a unique $t \Mod{p} : f(r + tp^j) \equiv 0 \Mod{p^{j + 1}}$.
\end{theorem}

\begin{definition}\label{polyx}
A \textit{Goldbach Polynomial} (G.P.) exists when it is possible to take some $a \in \mathbb{N}$ and for each $p_i \in \mathcal{P}_a$ a unique $\alpha_i \in \mathbb{N} \cup \{0\}$ exists producing a mapping $\mathcal{G}_- : \mathbb{C} \to \mathbb{C}$ where  
\begin{equation}\label{roots}
 \mathcal{G}_-(z) = \prod_{p_k \in \mathcal{P}_a} (z - p_k) - \prod_{p_k \in \mathcal{P}_a}p_k^{\alpha_k}
 \end{equation}
 and $\mathcal{G}_-(2a) = 0$. The roots are given by $\mathcal{R}_a = \{ r_m \in \mathbb{C} : \mathcal{G}_-(r_m) = 0 \}$. Note, for any $p_i \in \mathcal{P}_a$
 \begin{equation}\label{p}
 \mathcal{G}_-(p_i) = \prod_{p_k \in \mathcal{P}_a} (p_i - p_k) - \prod_{p_k \in \mathcal{P}_a}p_k^{\alpha_k} = (p_i -p_i)\prod_{p_k \in \mathcal{P}_a \setminus \{p_i\}} (p_i - p_k) - \prod_{p_k \in \mathcal{P}_a}p_k^{\alpha_k} = - \prod_{p_k \in \mathcal{P}_a}p_k^{\alpha_k}.
 \end{equation}
 \end{definition}

\begin{remark}\label{pderiv}
The Derivative of \ref{roots} follows via the product rule
\begin{equation}\label{deriv}
 \mathcal{G}_-'(z) = \prod_{p_k \in \mathcal{P}_a \setminus \{p_1\}} (z - p_k)  +\cdots + \prod_{p_k \in \mathcal{P}_a \setminus \{p_{\pi(a)}\}} (z - p_k)    = \sum_{p_s \in \mathcal{P}_a} \prod_{p_k \in \mathcal{P}_a \setminus \{p_s\}} (z - p_k). 
 \end{equation}
 Similarly to equation \ref{p}, for any $p_i \in \mathcal{P}_a$ at $z = p_i$ only one term in equation \ref{deriv} remains non-zero since $p_i$ is only absent from one of the product terms above, showing
 \begin{equation}\label{roots1}
  \mathcal{G}_-'(p_i) =  \prod_{p_k \in \mathcal{P}_a \setminus \{p_i\}} (p_i - p_k)\quad \text{where} \quad  \mathcal{G}_-'(p_i)  \equiv \mathcal{G}_-'(0)  \prod_{p_k \in \mathcal{P}_a \setminus \{p_i\}} (- p_k)\not \equiv 0 \Mod{p_i}.
 \end{equation} 
 Since $2a - p_k > 1$ for each $p_k \in \mathcal{P}_a$ shows $\mathcal{G}_-'(2a) > 1$ in \ref{deriv}, Theorem \ref{repeated} states $2a$ is not a repeated root and there exists a $Q(z) \in \mathbb{Z}[z]: \mathcal{G}_-(z) = (z - 2a)Q(z)$ where $Q(2a) \neq 0$ and $\mathcal{G}_-(p_i) = (p_i - 2a)Q(p_i)$. In Definition \ref{polyx}, $\mathcal{G}_-(2a) = 0$ shows $\prod_{p_k \in \mathcal{P}_a} (2a - p_k) = \prod_{p_k \in \mathcal{P}_a}p_k^{\alpha_k}$ where $(p_i - 2a) \prod_{p_k \in \mathcal{P}_a \setminus \{p_i\}} (2a - p_k) = - \prod_{p_k \in \mathcal{P}_a}p_k^{\alpha_k}$. Equation \ref{p} then allows via transitivity
 \begin{equation}\label{factor}
\mathcal{G}_-(p_i) = - \prod_{p_k \in \mathcal{P}_a}p_k^{\alpha_k} = (p_i - 2a)Q(p_i) = (p_i - 2a)\prod_{p_k \in \mathcal{P}_a \setminus \{p_i\}}(2a - p_k). 
\end{equation}
\end{remark}

\begin{example}
Let $a = 3 : \mathcal{P}_3 = \{2, 3\}$ and $\mathcal{G}_-(z) = (z-2)(z-3) - 2^2 \times 3$ with $\mathcal{R}_3 = \{6, -1\}$. Note, that the root of interest is of the form $6 = 2^2 + 2 = 3+3$. 
\end{example} 
 
 \subsection{A Proof of the Goldbach Conjecture}

\begin{proposition}\label{gbf}
\textit{Let $a \in \mathbb{N}_{> 3}$. If $2a$ is a counterexample to the G.C. then there exists a G.P.}
\end{proposition}

\begin{proof}
Assume $2a$ is counter-example to the G.C. Then no $2a - p_k$ may be prime when $p_k \in \mathcal{P}_a$. If any $2a - p_k$ had a prime factor $a < p < 2a$ then $a < 2a - p_k = n p < 2a$ with $p > a$ forcing $0< n < 2$. Hence, $2a - p_k = p$, contradicting the assumption that $2a$ is a counter-example to the G.C. Therefore, any counter-example requires each $2a - p_k$ only to be divisible by primes in $\mathcal{P}_a$ where the Fundamental Theorem of Arithmetic ensures unique $\alpha_1, \dots, \alpha_{\pi(a)} \in \mathbb{N}\cup \{0\}$ exist where $\prod_{p_k \in \mathcal{P}_a} (2a - p_k) = \prod_{p_k \in \mathcal{P}_a}p_k^{\alpha_k}$ and these $\alpha_i$ produce a G.P. from Definition \ref{polyx}.
\end{proof}

 We now proceed under the assumption that $a \in \mathbb{N}_{>3}$ and fixed where $2a$ a counter-example to the G.C. We then prove no G.P.s exist for $a > 3$ proving the G.C. is true via Proposition \ref{gbf}.
 
 \begin{lemma}\label{iff}
\textit{$\alpha_i > 0$ if and only if there exists at least one $p_m \in \mathcal{P}_a : 2a \equiv p_m \Mod{p_i}$.}
\end{lemma}

\begin{proof}
Since $2a$ is a root of the G.P. in Definition \ref{polyx}, it must satisfy  $\prod_{p_k \in \mathcal{P}_a} (2a - p_k) = \prod_{p_k \in \mathcal{P}_a}p_k^{\alpha_k}$ where $\alpha_i > 0$ iff $p_i$ divides at least one $2a - p_m$ term which is true iff $2a \equiv p_m \Mod{p_i}$. 
\end{proof}

\begin{definition}\label{pr}
Assume $p_i \in \mathcal{P}_a$. Equation \ref{factor} shows $\mathcal{G}_-(p_i) = - \prod_{p_k \in \mathcal{P}_a}p_k^{\alpha_k}$ where it follows $v_{p_i}(\mathcal{G}_-(p_i)) = \alpha_i$ and $\mathcal{G}_-(p_i) \equiv 0 \Mod{p_i}$ iff $\alpha_i > 0$. Under equation \ref{roots1} of Remark \ref{pderiv} it was shown $\mathcal{G}_-'(p_i) \not \equiv 0 \Mod{p_i}$. Therefore, via Hensel's Lemma \ref{HL}, if $\alpha_i > 0$,
\begin{equation}\label{pr0}
\mathcal{G}_-(p_i) \equiv 0 \Mod{p_i^{\alpha_i}} \quad \text{and} \quad \mathcal{G}_-(p_i + t_{\alpha_i} p_i^{\alpha_i}) \equiv 0 \Mod{p_i^{\alpha_i + 1}}
\end{equation} 
for a unique $p_i$-adic root $r_i \equiv p_i + t_{\alpha_i}p_i^{\alpha_i} \Mod{p_i^{\alpha_i + 1}}$ and is the \textit{Principal Root Modulo $p_i$} (P.R.). 
\end{definition}

\begin{remark}\label{func}
Definition \ref{pr} states there exists a unique P.R. satisfying equation \ref{pr0} if $\alpha_i > 0$. From Remark \ref{pderiv} it was shown $\mathcal{G}_-(z) = (z - 2a)Q(z)$ that equation \ref{pr0} in Definition \ref{pr} allows 
\begin{equation}\label{series1}
\mathcal{G}_-(p_i) \equiv (p_i - 2a)Q(p_i) \equiv 0 \Mod{p_i^{\alpha_i}}.
\end{equation}
along with
\begin{equation}\label{series2}
\mathcal{G}_-(p_i + t_{\alpha_i}p_i^{\alpha_i}) \equiv (p_i + t_{\alpha_i}p_i^{\alpha_i}-2a)Q(p_i + t_{\alpha_i}p_i^{\alpha_i}) \equiv 0 \Mod{p_i^{\alpha_i +1}}.
\end{equation}
From Remark \ref{pderiv} a cancellation of $p_i - 2a$ in equation \ref{factor} gives $Q(p_i) = \prod_{p_k \in \mathcal{P}_a \setminus \{p_i\}}(2a - p_k).
$ Using Definition \ref{polyx} and $Q(2a) \neq 0$ from Remark \ref{pderiv}, $Q(z)$ may be written in terms of its roots as
\begin{equation}\label{Q}
Q(z) = \prod_{r_k \in \mathcal{R}_a \setminus \{2a\}}(z - r_k)^{\sigma_k} \quad \text{where} \quad Q(p_i) = \prod_{r_k \in \mathcal{R}_a \setminus \{2a\}}(p_i - r_k)^{\sigma_k} = \prod_{p_k \in \mathcal{P}_a \setminus \{p_i\}}(2a - p_k).
\end{equation}
This will aid in proving $2a$ is always the P.R. when $\alpha_i > 0$ in the following lemma.  
\end{remark}

\begin{lemma}\label{perfect}
\textit{For any $p_i \in \mathcal{P}_a$ it follows that if $\alpha_i > 0$, then $2a$ is the Principal Root.}
\end{lemma}

\begin{proof}
Assume $\alpha_i > 0$ and, for contradiction, $v_{p_i}(2a - p_i) = 0$. This requires $Q(p_i) \equiv 0 \Mod{p_i^{\alpha_i}}$ from equations  \ref{factor} and \ref{pr0} when $\alpha_i > 0$. From equation \ref{series2} the  $v_{p_i}(p_i + t_{\alpha_i}p_i^{\alpha_i}-2a) = 0$ where $Q(p_i + t_{\alpha_i}p_i^{\alpha_i}) \equiv 0 \Mod{p_i^{\alpha_i +1}}$ for a unique P.R. in Definition \ref{pr}. The concern is not with the construction of the P.R., but instead with the $p_i$-adic representation of $2a$. Under Remark \ref{pderiv}, $Q(2a) \neq 0$, however, with $v_{p_i}(2a - p_i) = 0$, Lemma \ref{iff} states there exists at least one $p_m \in \mathcal{P}_a \setminus \{p_i\}$ where $2a \equiv p_m \Mod{p_i}$ making $2a$ a $p_i$-adic root of equation \ref{Q}. Since $2a$ is a global root and integer, $2a$ is contained in $Q(z)$ when $\alpha_i > 0$, contradicting $Q(2a) \neq 0$. Therefore, equation \ref{series2} requires $2a \equiv p_i + t_{\alpha_i}p_i^{\alpha_i} \Mod{p_i^{\alpha_i + 1}}$ to ensure consistency of Lemma \ref{iff} and equation \ref{series1}.
\end{proof}

\begin{theorem}\label{gbt}
\textit{The Goldbach Conjecture is true.}
\end{theorem}

\begin{proof}
Under Lemma \ref{perfect} if $\alpha_i > 0$, then Definition \ref{pr} shows $2a \equiv p_i + t_{\alpha_i}p_i^{\alpha_i} \Mod{p_i^{\alpha_i + 1}}$. Since the $p_i$-adic valuation of $2a - p_i$ is $\alpha_i$ whenever $\alpha_i > 0$, under equation \ref{factor} no prime may divide multiple $2a - p_k$ when $\alpha_i > 0$. Since all $2a - p_k > 1$ and must be a composition of only primes in $\mathcal{P}_a$, the pigeonhole principal shows for each $p_i \in \mathcal{P}_a$ that $2a = p_i^{\alpha_i} + p_i$. Furthermore, with each $2a- p_i > 1$, each $\alpha_i > 0$ and $a\# | 2a$. This is false when $a > 3$ since $2a < a\#$ for $a > 4$, proving no G.P.s exist beyond this value and the GC. is true via Proposition \ref{gbf}.
\end{proof}

\thispagestyle{plain}

\end{document}